\documentclass[11pt,reqno]{amsart}
\usepackage{amssymb, amsmath}
\usepackage{mathrsfs}
\usepackage{enumerate}

\newcommand {\PP}{{I\kern-.3em P}}

\newcommand {\HH}{{\mathbb{H}}}
\newcommand {\RR}{{\mathbb{R}}}

\newcommand{\beq}{\begin{equation}}
\newcommand{\eeq}{\end{equation}}
\newcommand{\beqq}{\begin{equation*}}
\newcommand{\eeqq}{\end{equation*}}

\numberwithin{equation}{section}

{\begin{enumerate}[\bfseries 1.]}{\end{enumerate}}
{\begin{enumerate}[\bfseries a.]}{\end{enumerate}}

\def\d{\displaystyle}

\newtheorem{thm}{Theorem}

\newtheorem{prop}{Proposition}
\newtheorem{cor}{Corollary}
\theoremstyle{remark}

\begin{document}
\title{Geodesic excursions into an embedded disc on a hyperbolic riemann surface}
\author{Andrew Haas }
 \email{haas@math.uconn.edu}
 \address{University of Connecticut, Department of Mathematics,
 Storrs, CT 06269}
 
\begin{abstract}

We calculate the asymptotic average rate at which a generic geodesic on a finite area hyperbolic 2-orbifold returns to an embedded disc on the surface, as well as the average amount of time it spends in the disc during each visit. This includes the case where the center of the disc is a cone point.
 \end{abstract}

 
\maketitle
 \markboth {Andrew Haas}{Geodesic excursions}

\section{introduction}

A finite area hyperbolic 2-orbifold $S$ is the quotient of the Poincar\'e upper-half plane $\HH$ by a discrete group $\Gamma\subset \text{PSL}_2(\RR).$ Let $\pi:\HH\rightarrow
S=\HH/\Gamma$ denote the covering projection and let $q\in S$ be a non-cone point of $S$. Without loss of generality, suppose that $i\in \HH$ is a preimage of $q$. As usual, $B_r(x)$ is the open disc of radius $r$ centered at $x$.  There is a number $R>0$  for which the map $\pi: B_R(i)\rightarrow B_R(q)$ is an embedding.

Each vector $v$ in the unit tangent bundle $T_1S$ uniquely determines a geodesic ray $\gamma=\gamma(v)$ which is the projection of the orbit of the geodesic flow, $\{G_t(v)\,|\,t\in [0,\infty)\},$ to $S$. If the ray keeps returning to the disc $B_r(q)$, then there is a sequence of pairs of parameter values $\{(t_i,s_i)\}$ so that 
\beqq
\bigcup_{i=1}^{\infty}\gamma([t_i,s_i])=\gamma\cap \overline{B}_r(q),
\eeqq
in other words, $\gamma$ meets the closed disc $B_r(q)$ in precisely the arcs $\gamma([t_i,s_i])$ of $\gamma.$ We shall refer to these arcs as the excursions of $\gamma$ into $B_r(q).$
For $v\in T_1S$ define the counting function $N_v(r)(t)=\max\{i\,|\,t_i\leq t\}.$ Then we have

\begin{thm}\label{return0}
For almost all $v\in T_1S$ and all $0<r<R$
\beq\label{average}
\lim_{t\rightarrow\infty}\frac{1}{t}N_v(r)(t)=\frac{2\sinh r}{ area (S)}.
\eeq
\end{thm}

Two  consequences of the theorem are given in

\begin{cor}\label{return1}
For almost all $v\in T_1S$ and all $0<r<R$
\begin{itemize}

\item[{\rm (i)}] The asymptotic average  length of an excursion is
\beq\label{part1}
\lim_{n\rightarrow\infty}\frac{1}{n}\sum_{i=1}^{n} (s_i-t_i)=\frac{\pi(\cosh r-1)}{\sinh r}.
\eeq

\item[{\rm (ii)}]
The average length of an arc joining the beginnings of consecutive excursions is 
\beq\label{part2}
\lim_{n\rightarrow\infty}\frac{1}{n}\sum_{i=1}^{n-1} (t_{i+1}-t_{i})=\frac{ area (S)}{2\sinh r}.
\eeq

\end{itemize}
\end{cor}

Note that when $r$ is large, (\ref{part1})   is close to $\pi.$ 
 In Section \ref{3} we expand these results to cover the case where $q$ is a cone point. We will also  consider their expressions  in terms of the  area of a disc rather than  its radius.

 \section{Details}

The unit tangent bundle of $\HH$ is $T_1\HH=\HH\times S^1.$ In these coordinates the geodesic flow $\tilde{G}$ has the natural invariant measure $\tilde{m}=dAd\theta.$ The non-vertical vectors in $T_1\HH$ are a set of full measure and can be parameterized by the set of triples $(x,y,t)\in\RR^3$ with $x\not=y.$ Let $\alpha$ be the geodesic in $\HH$ oriented from the  endpoint $\alpha_{-}=y$ to $\alpha_{+}=x $  and parameterized so that $\alpha(0)$ is the Euclidean midpoint of the semicircle $\alpha(\RR)$. Then $(x,y,t)$ corresponds to the vector $v=\dot{\alpha}(t)\in T_1\HH.$ In these coordinates the geodesic flow has the form $\tilde{G}_s(x,y,t)=(x,y,t+s).$ Furthermore, the geodesic flow on $T_1S$ has the invariant probability measure $\mu$, whose lift to $T_1\HH$ is equal to
\beq
\tilde{\mu}= \frac{1}{2\pi\,\text{area}(S)}\tilde{m}=\frac{1}{\pi\, \text{area}(S)}\varphi(x,y)\, dx dy dt.
\eeq
where $\varphi (x,y)=(x-y)^{-2}$ \cite{nicholls}.  

Let $\mathcal{E}\subset T_1S$ be the set of vectors $v$ for which there are infinitely many excursions of $\gamma=\gamma(v)$ into $B_r(q)$ for each $0<r<R.$ As a consequence of  the Poincar\'e recurrence Theorem,  \cite{sinai}, $\mathcal{E}$ has full measure in $T_1S$. Observe that if $v\in \mathcal{E},$ then $G_t(v)\in \mathcal{E}$ for all $t\in \RR.$

Let $\mathcal{L}_r^*$ be the subset of the unit tangent bundle over the circle $\partial B_r(q)$,  (bounding the disc $B_r(q)$)  consisting of vectors that point into the disc $B_r(q).$  Then $\mathcal{L}_r=\mathcal{L}_r^*\cap\mathcal{E}$ is a cross-section for the geodesic flow on $T_1S$, \cite{bks}. In other words, for almost all $v\in T_1S$ there exists an increasing sequence of values $\tau_i$ so that $G_{\tau_i}(v)\in \mathcal{L}_r$.
Given $\epsilon<R-r$, define the $\epsilon$-thickened section $ \mathcal{L}_r(\epsilon)=\{G_t(v)\,|\, t\in [0,\epsilon], v\in \mathcal{L}_r\}.$ Analysis of the thickened section is the main tool in the proof of Theorem \ref{return0}. This is similar to the approaches taken in \cite{haas}, \cite{ nakada} and \cite{strat}. 

\begin{prop}\label{thick}
\beq
\mu( \mathcal{L}_r(\epsilon))=\frac{2\epsilon\sinh r}{ area (S)}
\eeq
\end{prop}

The cross-section $\mathcal{L}_r$ lifts to a subset $\tilde{\mathcal{L}}_r$ of $T_1\HH$ over the circle $\partial B_r(i)$ in $\HH$. Then $ \mathcal{L}_r(\epsilon)$ has a connected preimage $\Lambda_r(\epsilon)\supset \tilde{\mathcal{L}}_r(\epsilon)$ in $T_1\HH.$ Furthermore, since $\pi$ is an embedding of $B_r(q)$ into $S$ and $\epsilon<R-r$, the projection $\pi_{*}:T_1\HH\rightarrow T_1S$ restricts to an embedding  of $\Lambda_r(\epsilon)$ onto $ \mathcal{L}_r(\epsilon)$. It follows that $\Lambda_r(\epsilon)=  \tilde{\mathcal{L}}_r(\epsilon).$ Thus $\mu( \mathcal{L}_r(\epsilon))=\tilde{\mu}( \tilde{\mathcal{L}}_r(\epsilon))$ and we can perform computations in $\RR^3$ using the coordinates given earlier.

Define the M\"obius transformations:
\beqq
W_{\rho}(x)=\d{\frac{z\sinh\rho -1}{z+\sinh\rho}}\,\,\,\,\,\,\,\text{and}\,\,\,\,\,\,\,  
U_{\rho}(x)=\d{\frac{z\sinh\rho +1}{-z+\sinh\rho}}.
\eeqq

For real numbers $a\not= b$ let $\overline{ab}$ denote the geodesic $\alpha$ in $\HH$ with endpoints $\alpha_+=a$ and $\alpha_-=b.$

We will need the following extension of Theorem 5 from \cite{haas3}.

\begin{thm}\label{coneold}
Given $x\geq 0$, the geodesics $\overline{xw}$ and $\overline{xu}$ with $w=W_{\rho}(x)\in [-1/x,x)$ and $u=U_{\rho}(x)$, are both tangent to the disc $B_{\rho}(i)$. Furthermore, given $x\leq 0$, the geodesics $\overline{xw}$ and $\overline{xu}$ with $w=-W_{\rho}(-x)\in (x,-1/x]$ and $u=-U_{\rho}(-x)$, are both tangent to the disc $B_{\rho}(i).$  

\end{thm}

Although the proof in \cite{haas3} only addressed the tangent geodesic with endpoints $x$ and $W_{\rho}(x)$, the other case can be uncovered there and follows by choosing $\xi=z(c+i\sqrt{1-c^2})$ in Lemma 3 and $\eta=c+i\sqrt{1-c^2}$ in the proof of Theorem 5.
\vskip .1in
  
\noindent {\em Proof of Proposition \ref{thick}}.\,
Fix $0<r<R$. For a given $x\in \hat{\RR}$, let $I_x$ denote the set of values $y\in \hat{\RR}$ for which the geodesic $\overline{xy}$ intersects the closed disc $\overline{B}_r(i)$. It follows from Theorem \ref{coneold} that for $x\geq 0,\, I_x$ is the interval between the points $W_{\rho}(x)$   and $U_{\rho}(x)$; the unbounded interval if $x \leq \sinh \frac{r}{2}$ and the bounded interval if  $x>\sinh  \frac{r}{2}$. A similar statement is true for negative $x$.

For $x\in \RR$ and $y\in I_x$, let $t_{xy}$ denote the parameter value for which the unit tangent vector   $(x,y,t_{xy})\in \tilde{\mathcal{L}}_r.$
Then
\beqq
 \tilde{\mathcal{L}}_r(\epsilon) =\{ (x,y,t)\,|\,x\in\RR,\, y\in I_x\, \text{and}\, t_{xy}\leq t\leq t_{xy}+\epsilon\}.
\eeqq
Consequently, $ \tilde{\mu}(\tilde{\mathcal{L}}_r(\epsilon))=$
\beq\label{area}
\frac{1}{\pi\,\text{area}(S)}\int_{\RR}\int_{I_x}\int_{t_{xy}}^{t_{xy}+\epsilon}\varphi(x,y)\,dtdydx = \frac{\epsilon}{\pi\,\text{area}(S)}\int_{-\infty}^{\infty}\int_{I_x}\varphi(x,y)\,dydx.  
\eeq
This unpacks as the following sum of integrals:
\beqq 
\frac{\epsilon}{\pi\,\text{area}(S)}\left [\int_0^{\sinh  \frac{r}{2}}\int_{-\infty}^{W_r(x)}\frac{1}{(x-y)^2}\,dydx+
\int_0^{\sinh  \frac{r}{2}}\int_{U_r(x)}^{\infty}\frac{1}{(x-y)^2}\,dydx+\right .
\eeqq

\beqq
\hskip .8in \int_{\sinh  \frac{r}{2}}^{\infty}\int_{U_r(x)}^{W_r(x)}\frac{1}{(x-y)^2}\,dydx+
\int_{-\sinh  \frac{r}{2}}^0\int_{-W_r(-x)}^{\infty}\frac{1}{(x-y)^2}\,dydx+
\eeqq
\beqq
\hskip 1in\left . \int_{-\sinh  \frac{r}{2}}^0\int_{-\infty}^{-U_r(-x)}\frac{1}{(x-y)^2}\,dydx+
\int_{-\infty}^{-\sinh  \frac{r}{2}}\int_{-W_r(-x)}^{-U_r(-x)}\frac{1}{(x-y)^2}\,dydx \right ]
\eeqq
The  integrals are easily computed.\\  
\qed

\noindent {\em Proof of Theorem \ref{return0}}.\,
 Let $\chi_Y(t)$ denote the characteristic function of the set $Y$. Then for $\epsilon<R-r$ we have the inequalities
 \beqq
 \int_0^t\chi_{\mathcal{L}_r(\epsilon)}(G_{\tau}(v))d\tau\,-2\epsilon\,\leq\, \epsilon N_v(r)(t)\,\leq\,
  \int_0^t\chi_{\mathcal{L}_r(\epsilon)}(G_{\tau}(v))d\tau\,+2\epsilon.
  \eeqq

Divide thru by $\epsilon t$. Since  the geodesic flow is ergodic, as $t\rightarrow\infty$ the left and right hand limits approach 
$(1/\epsilon)\mu(\mathcal{L}_r(\epsilon)).$
The theorem then follows from Proposition \ref{thick}.
\qed

\vskip .1in

\noindent {\em Proof of Corollary \ref{return1}}.\,
First we  address part (ii). Let $v\in\mathcal{E}$ and let $\{(t_i,s_i)\}$ be the sequence of excursions into the disc $B_r(q)$ associated to the geodesic $\gamma=\gamma(v).$ In particular we have $N_v(r)(t_n)=n$ (or $n-1$ if $\gamma(0)\in B_r(q))$. Then by Theorem \ref{return0},  for almost all $v\in T_1S$
\beqq\label{coravg}
\lim_{n\rightarrow\infty}
\frac{1}{n}\sum_{i=1}^{n-1}(t_{i+1}-t_i)=
\lim_{n\rightarrow\infty}\frac{t_n}{n}=
\lim_{n\rightarrow\infty}\frac{t_n}{N_v(r)(t_n)}
=\frac{\text{area}(S)}{2\sinh r},
\eeqq
 proving part (ii).

For part (i) observe that the limit (\ref{part1}) can be expressed as an integral:
\beq\label{flow}
\lim_{n\rightarrow\infty}\frac{1}{n}\sum_{i=1}^{n} (s_i-t_i)=\lim_{n\rightarrow\infty}\frac{1}{n}\int_0^{s_n}\chi_{(B_r(q))}\left(G_v(\tau)\right)d\tau
\eeq
By Theorem \ref{return0} and the Ergodic Theorem for flows, this can be further divided into a product of limits, each piece of which converges for almost all $v\in T_1S$. Thus (\ref{flow}) is equal to 
\beq\label{split}
\begin{split}
& \lim_{n\rightarrow\infty}\frac{1}{s_n}\int_0^{s_n}\chi_{(B_r(q))}\left(G_v(\tau)\right)d\tau\,\times\, 
\lim_{n\rightarrow\infty}\frac{s_n}{n}\\
\\
& =\frac{ \text{area}(B_r(q))}{\text{area}(S)}\,\times\,\frac{\text{area}\,(S)}{2\sinh r}=\frac{\pi(\cosh r-1)}{\sinh r}
\end{split}
\eeq
as asserted in part (i) of the Corollary, \cite{Beardon}. 
\qed

\section{Cone points}\label{3}
 
If the stabilizer of $i$ in $\Gamma$ is a subgroup of order $k>0$ then $q$ is called a cone point of order $k$. Let $R>0$ be the largest value so that for all $0<r<R$ the projection  $\pi: B_r(i)\rightarrow B_r(q)$ is precisely $k-\text{to}-1$ in the complement of $i$. Then for $0<r<R$, the definition of $N_v(r)(t)$ extends without modification to the case where $q$ is a cone point. The thickened sectioned is defined as before but now,  since  the map $\pi_*$ is $k-\text{to}-1$, instead of Proposition \ref{thick}, we have
\beq\label{withcone}
\frac{1}{\epsilon}\mu( \mathcal{L}_r(\epsilon)= \frac{1}{k\epsilon}\tilde{\mu}( \tilde{\mathcal{L}}_r(\epsilon))=\frac{2\sinh r}{k\,\text{area}(S)}.
\eeq
Thus Theorem \ref{return0} remains valid when $q$ is a cone point of order $k$, with the right hand side of (\ref{average}) replaced by (\ref{withcone}).

When  studying the rate of geodesic return to the neighborhood of a cone point, it is  interesting to have all the quantities   expressed  in terms of area rather than radius. Let $D_a(q)$ denote the disc of area $a$ about $q$. 
Observe that, if $q$ is a cone point of order $k$, then the area of the disc $B_r(q)$  is 
\beqq
\text{area}(B_r(q))=\frac{1}{k}\text{area}(B_r(i))= \frac{2\pi}{k} (\cosh r -1). 
\eeqq
Thus $D_a(q)$ is a disc or radius $r=\cosh^{-1}(\frac{k}{2\pi}a+1).$ 

Let  $\{(t_i,s_i)\}$ be the sequence of excursions of the geodesic $\gamma=\gamma(v)$ into the disc $D_a(q) $. Then, almost as before, for $v\in\mathcal{E}$ we can define $\mathcal{N}_v(a)(t)=\max\{i\,|\,t_i\leq t\}.$ In this setting  Theorem \ref{return0} for a cone point of order $k$,  expressed in terms of area, takes the form:
\begin{thm}\label{area}
For almost all $v\in T_1S$ and all $0<a<A=(2\pi/k) (\cosh R -1)$
\beq\label{average2}
\lim_{t\rightarrow\infty}\frac{1}{t}\mathcal{N}_v(a)(t)=\frac{2\sqrt{(\frac{1}{k}+\frac{1}{2\pi}a)^2+\frac{1}{k^2}}}{\text{area}\,(S)}.
\eeq
\end{thm}
When $k$ is large this value is close to  $a/(\pi\,\text{area}(S))$, which  it the value computed in \cite{haas} for the asymptotic rate of return to a cusp neighborhood.

Also,   the asymptotic average length of an excursion can still be computed when $q$ is a cone point. Replacing $B_r(q)$ by $D_a(q)$, equation (\ref{flow}) holds by definition and, as in equation (\ref{split}), using Theorem \ref{area} and  the Ergodic Theorem we have
\begin{cor}\label{areaavg}
For almost all $v\in T_1S$ and all $0<a<A$, the average length of an excursion into $D_a(q)$ is
\beqq
\lim_{n\rightarrow\infty}\frac{1}{n}\sum_{i=1}^{n} (s_i-t_i)=  \frac{a}{\text{area}(S)}\times
\frac{\text{area}(S)}{2\sqrt{(\frac{1}{k}+\frac{1}{2\pi}a)^2+\frac{1}{k^2}}}=
\frac{\pi}{\sqrt{(\frac{2\pi}{ka}+1)^2+(\frac{2\pi}{ka})^2}}
\eeqq
\end{cor}

For $k$ large the average length of an excursion is close to $\pi$, which was the value  computed in \cite{haas} for the average length of an excursion relative to a cusp neighborhood. When the area of the disc is large, the average length of an excursion is also close to $\pi,$  independent of the order of the cone point.

\end{document}